\documentclass[10pt,final,onefignum]{siamltex}
\usepackage{amsmath,pstricks}
\usepackage{hyperref}
\usepackage{bbm}
\usepackage{graphicx}

\newcommand{\Div}{\nabla\cdot}

\newcommand{\ep}{\varepsilon}
\newcommand{\im}{\mathtt{i}}
\newcommand{\e}{\mathtt{e}}
\newcommand{\dO}{\partial\Omega}
\newcommand{\D}{\mathbbm{D}}
\newcommand{\T}{\mathbbm{T}}
\newcommand{\R}{\mathbbm{R}}
\newcommand{\N}{\mathbbm{N}}
\newcommand{\X}{\mathbbm{X}}
\newcommand{\tZ}{\widetilde{Z}_n^N}
\newcommand{\hZ}{\widehat{Z}_n^N}
\newcommand{\tv}{\underline{\mathrm{\tau}}}
\newcommand{\nv}{\underline{\mathrm{n}}}
\newcommand{\ud}{u_\delta}
\newcommand{\zd}{z_\delta}
\newcommand{\Ud}{U_\delta}
\newcommand{\od}{w_\delta}
\newcommand{\Gd}{G_\delta}
\newtheorem{remark}[theorem]{Remark}
\title{On the existence and uniqueness of weak solutions for 
  a vorticity seeding model}  
\author{Luigi C.~Berselli\thanks{Dipartimento di Matematica
    Applicata ``U.~Dini,'' Universit\`a di Pisa, Via B. Pisano 25/b
    I-56126 Pisa, Italia. Phone: +39 050 2219479, Fax: +39 050
    2219451, email: {\tt berselli@dma.unipi.it}}
        \and 
Marco Romito\thanks{Dipartimento di Matematica, Universit\`a di
    Firenze, Viale Morgagni 67/a I-50134 Firenze, Italia  
Phone: +39 055 4237136, Fax +39 055 4222695, email: {\tt
    romito@math.unifi.it}}} 
\begin{document}
\maketitle
\begin{abstract}
In this paper we study the Navier-Stokes equations with a Navier-type
boundary condition that has been proposed as an alternative to common
near wall models. The boundary condition we study, involving a linear
relation between the tangential part of the velocity and the
tangential part of the Cauchy stress-vector, is related to the 
vorticity seeding model introduced in the computational approach to 
turbulent flows. The presence of a point-wise non vanishing normal
flux may be considered as a tool to avoid the use of phenomenological 
near wall models, in the boundary layer region. Furthermore, the
analysis of the problem is suggested by recent advances in the study
of Large Eddy Simulation.

In the two dimensional case we prove existence and uniqueness of weak
solutions, by using rather elementary tools, hopefully understandable 
also by applied people working on turbulent flows. The asymptotic
behaviour of the solution, with respect to the averaging radius
$\delta,$ is also studied. In particular, we prove convergence of the
solutions toward the corresponding solutions of the Navier-Stokes
equations with the usual no-slip boundary conditions, as the small
parameter $\delta$ goes to zero.
\end{abstract}
\pagestyle{myheadings}
\thispagestyle{plain}
\markboth{LUIGI C. BERSELLI AND MARCO ROMITO}{ON A VORTICITY SEEDING MODEL} 
\section{Introduction}
In this paper we consider the Navier-Stokes equations and in
particular the role of
boundary conditions in the simulation of boundary effects in
turbulent flows. We consider the Navier--Stokes equations (in
non-dimensional form) for viscous incompressible fluids in a bounded
smooth domain $\Omega\subset\R^n,$ $n=2,3$, 
\begin{equation}\label{main}
\left\{\begin{array}{ll}
\partial_t u-\displaystyle\frac{2}{Re}\Div\D(u)+(u\cdot\nabla)u+
\nabla p=0\quad\text{in }\Omega\times(0,T)\\\\
\Div u=0\quad\text{in }\Omega\times(0,T).\\
\end{array}\right.
\end{equation}
We recall that $u=(u_1,\dots,u_n)$ is the unknown velocity field, $p$ 
is the hydrostatic pressure, $Re>0$ is the Reynolds number, and
$\D(u)$ is the deformation tensor, \textit{i.e.}, the symmetric part
of the matrix of derivatives of $u$:
$$
\D(u)=\frac12\left(\frac{\partial u_i}{\partial
    x_k}+\frac{\partial u_k}{\partial x_i}\right),
$$ 
and the Navier-Stokes equations are generally equipped with the
\textit{no-slip} boundary conditions on $\Gamma=\dO$ 
\begin{equation}\label{Dirichlet}
  u=0\quad\text{on }\Gamma\times(0,T).
\end{equation}
To introduce the problem that we will study, we recall that while at a
free surface it is natural to require continuity of the stress-tensor
($\mathbbm{I}$ being the identity) 
$$
\T(u,p)=-p\,\mathbbm{I}+\frac{2}{Re}\D(u),
$$
the conditions at a solid wall are much troublesome. The no-slip
condition~\eqref{Dirichlet} has been justified by
Stokes~\cite{Stok1845} since the contrary assumption 
\begin{quote}
\textrm{\dots implies an infinitely greater resistance to the sliding
  of one portion of fluid past another than to the sliding of fluid
  over a solid}.  
\end{quote}
 It is well-known that there are situation in which
the boundary condition~\eqref{Dirichlet} may not be valid. For
instance in Serrin~\cite{MR21:6836b} \S 64 it is pointed out that in
high altitude aerodynamics, or more generally when moderate pressure
and low surface stresses are involved, the adherence condition is no
longer true, see also the review in Truesdell~\cite{MR13:794d}. In
this respect several authors proposed various \textit{slip} (generally
nonlinear) conditions, modeling precise physical situations; see for
instance Serrin~\cite{MR21:6836b}, Beavers and Joseph~\cite{bj}, Krein
and Laptev~\cite{MR40:1714}.\par 
From the historical point of view, the slip (with friction) boundary
condition proposed\footnote{Note that in contrast to Stokes (1845)
  that used the continuum mechanics, Navier (1823) derived the
  equations by using some formal analogy with the elasticity theory
  and the assumption that molecules are animated by attractive and
  repulsive forces.} 
by Navier~\cite{Nav23} was     
\begin{equation}\label{navier}
u\cdot\nv=0\quad\text{and}\quad \beta\,
u_\tau+\underline{\cal T}(u,p)=0,\quad\beta>0,\quad\text{on
}\Gamma\times(0,T),      
\end{equation}
where $\nv$ denotes the unit normal vector to $\Gamma,$
$u_\tau=u-(u\cdot\nv)\nv,$ while $\underline{\cal
  T}(u,p)=\underline{t}(u,p)- 
(\underline{t}(u,p)\cdot\nv)\,\nv$ 
  denotes the tangential part of the \textit{Cauchy stress vector}
  $\underline{t}$ defined by: 
$$
\underline{t}(u,p)=\nv\cdot\T(u,p)=\sum_{k=1}^n \T_{ik}(u,p)\nv_k.
$$
The controversial between the condition~\eqref{navier} proposed by
Navier, versus the Dirichlet~\eqref{Dirichlet} proposed by Stokes was
analysed also by Maxwell, who observed that the same conditions may be 
derived in the kinetic theory of gases and that the parameter $\beta$
should depend on the Reynolds number $Re$ and on the mean free-path
$\lambda,$  satisfying the couple of consistency conditions
\begin{equation*}
  \begin{aligned}
    \beta&\to\infty\text{ as }\lambda\to0,\text{for $Re$ fixed}\\
    \beta&\to0\text{ as }Re\to\infty,\text{for $\lambda$ fixed}.
  \end{aligned}
\end{equation*}
With the above asymptotics it possible to recover in both cases the
correct no-slip boundary condition for viscous fluids and the
no-penetration condition for ideal fluids. A study of the numerical
problems related with the implementation of~\eqref{navier} can be
found in John~\cite{MR2003h:76033}. 

Recently, Fujita~\cite{MR2003m:76047} performed the analysis with
the ``slip or leak with friction'' boundary conditions. These
conditions are of particular interest in the study of polymers, blood
flow, and flow through filters. The boundary conditions studied
in~\cite{MR2003m:76047} with the techniques of variational
inequalities, turn out to be particular cases of the nonlinear
boundary condition proposed at p.~240 in~\cite{MR21:6836b} and they
are very strictly related to both the Navier and the no-slip boundary
conditions. See also Consiglieri~\cite{MR2004b:76006} for similar
problems.

Among other nonstandard boundary conditions we recall that studied
by Begue~\textit{et al.}~\cite{MR90j:35161} and the ``do-nothing''
Neumann conditions, very appealing for numerical studies,  implemented
in Heywood~\textit{et al.}~\cite{MR97f:76045}.
\subsection{Near Wall Models and turbulent flows}
Our interest in non standard boundary conditions comes essentially
from the study of turbulent flows.
First, we recall that in the boundary layer theory several log-law and
power-law are introduced, together with the fictitious boundaries, in
order to model turbulent flows within a small region near the
boundary. Roughly speaking, appropriate nonlinear boundary conditions
are imposed on an artificial boundary that lies inside the
computational domain. The boundary conditions may simulate (at least
in a computational approach) the behaviour of the boundary layer and
they are modeled to take into account the peculiar behaviour of a
fluid near the boundaries\footnote{In this respect we quote
  J.C. Maxwell: ``It is almost certain the the stratum of gas next to a
  solid body is in a very different state from the rest of the 
  gas''.}

Our main motivation comes from the mathematical theory of Large Eddy
Simulation (briefly, LES). In fact, the purpose of LES is to model
the evolution of large coherent structures (eddies); this is done by
studying the equations satisfied by a filtered velocity. Generally,
the filtered velocity $\overline{u}$ is defined through a convolution  
$$
\overline{u}(x,t)=g_\delta(x)*u(x,t)
$$
with a rapidly decreasing smoothing-kernel $g_\delta$ of width
$\delta;$  in several cases of practical interest 
$$
g_{\delta}(x)=\left(\displaystyle\frac{\gamma}{\pi}\right)^{3/2}
\displaystyle\frac{1}{\delta^3}
\ \textrm{e}^{\,-\textstyle\frac{\gamma|x|^2}{\delta^2}}
$$ 
By definition, the value of $\overline{u}$ at a
point $x_0$ on the boundary $\Gamma$ will depend on the behaviour of
$u$ in a neighborhood of width $\delta$ near that point; even if $u$
is extended to zero for each $x\not\in\overline{\Omega},$ it is clear
that in general $\overline{u}(x_0)\not=0$.

As pointed out in Galdi and Layton~\cite{MR2001h:76074} the physical
intuition may suggest that large coherent structures touching a wall
do not penetrate, sliding along it, and losing their energy. The
boundary condition of Navier may be revisited, by linking 
the micro-scale $\lambda$ of the kinetic theory of gases, with the
radius $\delta$ of the averaging filter.

Many Near Wall Models (NWM) have been tested in the computational
approach, see Sagaut~\cite{MR2003j:76070} and Piomelli and 
Balaras~\cite{MR2003b:76103}. The results are not uniformly successful 
and a positive application is very often based on a fine tuning of
parameters. This is why new models require at least a positive
background from the physical hypotheses and a coherent mathematical
analysis. In particular, a successful application of the
Navier slip-with-friction boundary condition~\eqref{navier} is
prevented by two main 
facts: 1) the presence of recirculation regions and 2) the presence of
fast time-fluctuating quantities.

The first problem is motivated by the fact that in recirculation
regions the local Reynolds number is very different from the main
stream and it is natural to expect that $\beta$ should depend
(possibly in a nonlinear way) on a local Reynolds number related to
the local slip speed, \textit{i.e.},  
$$
\beta=\beta(\delta,\,|u_\tau|).
$$
Preliminary analysis has been performed by John \textit{et
  al.}~\cite{jls04}, Dunca \textit{et al.}~\cite{djls2001}, and an
appropriate power-law choice of $\beta$ seems  promising to improve
the estimation of reattachment points. 
 
The limitation of Navier law~\eqref{navier} in a Boundary Layer theory
is that it can well describe time-averaged flow profiles, while the
information coming from fluctuating quantities in the wall-normal 
direction can play an important role in in triggering separation and
detachment.  
To try to overcome this limitation, very 
recently Layton~\cite{lay2002} recognized a particular class of
boundary conditions, leading to conditions similar ``in spirit'' to
the so called \textit{vorticity seeding methods}. In fact, a Navier
slip-with-friction boundary condition implies the generation of
vorticity at the boundary, proportional to the tangential
velocity. More precisely, in the case of a two-dimensional domain
$\Omega,$ for each smooth function $v$ such that  $v\cdot \nv=0$ on
the boundary, it holds  
$$
\nv\cdot\D(v)\cdot\tv-\frac{1}{2}\textrm{curl}\,v+k(v\cdot\tau)=0
\qquad\text{on }\Gamma, 
$$
where $\text{curl}\,v=\partial v_1/\partial x_2-\partial v_2/\partial
x_1,$ $\tv$ denotes the unit tangent vector, while $k$ is the
curvature of $\Gamma.$ 

In particular, in~\cite{lay2002} the following boundary condition
is proposed to simulate the boundary effects  
\begin{equation}\label{normal_component}
u\cdot\nv=\delta^2 g(x,t)\quad\text{and}\quad \frac{L}{\delta Re}
u_\tau+\underline{\cal T}(u,p)=0,  
\end{equation}
where $g$ is a \textit{highly oscillating} function in the time
variable (hopefully a random variable in numerical tests), while
it may be very smooth in the space variables and it should
satisfy the natural compatibility condition
\begin{equation}\label{compa}
\int_\Gamma g(x,t)\,d\sigma=0\qquad\forall\,t\in(0,T),
\end{equation}
which is required by  the normal trace of a divergence-free
vector field.

This way of reasoning is also similar to the introduction of
stochastic fluctuations, to simulate the micro-scale effects. A
comprehensive introduction to stochastic partial differential
equations in fluid mechanics and the statistical approach can be found
in Monin and Yaglom~\cite{MY81} and one main mathematical paradigm is
that an additional non-smooth term on the right-hand side may
naturally take into account the effect of the fast fluctuating
quantities. This leads to study the system 
$$
\partial_t u-\frac{2}{Re}\Div\D(u)+(u\cdot\nabla)u+\nabla
p=f+\frac{\partial g}{\partial t},
$$
where $g$ is a function that does  not have a proper time derivative,
but is just continuous or with other weak properties. Appropriate
notion of solutions, together with the statistical properties of the
above problem are studied in Bensoussan and Temam~\cite{MR50:1336} and
Vi{\v{s}}ik and Fursikov~\cite{MR82g:35095}. 
\subsection{Setting of the problem}
In the sequel we will restrict to the two dimensional case, since the
nonlinear character of the equations imposes some
restriction, see Remark~\ref{2D-3D}. Furthermore,
we fix the values of $L$ and of the Reynolds number to $2$ due to the
fact that we will not deal with the singular limit $Re\to\infty$
(regarding this limit see also the recent works of Clopeau~\textit{et
  al.}~\cite{MR99g:35102} and Mucha~\cite{MR1999576}). In our case we
will study the following boundary-initial value problem:
\begin{equation}\label{maina}
\begin{cases}
\partial_t u-\Div\D(u)+(u\cdot\nabla)u+\nabla p=f&\text{in
}\Omega\times(0,T)\\ 
\Div u=0&\text{in }\Omega\times(0,T)\\
u\cdot\nv=\delta^\alpha g(x,t)&\text{on }\Gamma\times(0,T)\\
u\cdot\tv+\delta\,\nv\cdot\D(u)\cdot\tv=0&\text{on
}\Gamma\times(0,T)\\
u(x,0)=u_0(x)&\text{in }\Omega.
\end{cases}
\end{equation}
Note that a similar problem, involving the
Smagorinsky-Lady\v{z}henskaya turbulence model together with a
nonlinear dependence of $\beta$ on $u,$ has been studied for instance 
by Par{\'e}s~\cite{MR95e:35160}, but in this reference the normal
datum $g$ is not allowed to depend on the time variable. In the sequel
our main interest will be to find weak hypotheses on $g(x,t)$ with
respect to the time variable (without any essential restriction on the
space regularity) that allow to prove existence of weak solutions to
the Navier-Stokes equations, see Theorem~\ref{mainex}. In particular,
in our analysis we will focus on two main points:
\begin{itemize}
\item[a)] to show the existence of weak solutions in the sense of
Leray-Hopf (since we do not want to deal with any weaker concept of
solution);
\item[b)] to use only elementary tools of functional analysis, in
  order to keep the paper intelligible to non specialists.
\end{itemize}
In other words, we want to consider solutions in an usual sense and we
want also to interact with applied people interested in this problem,
still keeping all the mathematical rigor needed to deal with the
problem and a certain sharpness of the results. In the case of
non-homogeneous no-slip conditions, several results of existence and
uniqueness of suitable solutions can be found in
Amann~\cite{MR2002b:76028}. 
 
Our intent to have a non-vanishing normal datum can be heuristically
understood also with the following argument: suppose that (for
simplicity in two dimensions) we have a fictitious boundary $\Gamma_1$
and we want to impose a condition on it in order to resolve
numerically the equation in a smaller domain $\Omega_1\subset\Omega,$
that rules out the boundary layer (see the figure below). 
\begin{figure}[h]
\centering
\begin{pspicture}(13,3.6)
\psframe[fillcolor=lightgray,fillstyle=solid](5,0.8)(8,2.4)
\rput(4.8,0.6){A}
\rput(8.1,0.6){B}
\rput(8.1,2.6){C}
\rput(4.8,2.6){D}
\rput(6.5,0.2){\Large $\Gamma$}
\rput(6.5,2.9){\Large $\Gamma_1$}
\rput(2.8,1.6){Boundary Layer}
\rput(10.2,1.6){Boundary Layer}
\pscurve(5,0.8)(2,1)(0.5,1.6)
\pscurve(8,0.8)(11,1.2)(12.5,1.7)
\pscurve(8,2.4)(9,2.6)(11,3.6)
\pscurve(5,2.4)(4,2.6)(3,3.6)
\end{pspicture}
\caption{The fictitious boundary}
\end{figure}
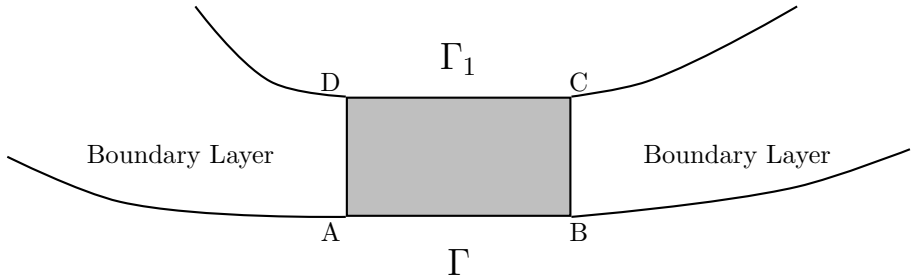

We have to require, by the incompressibility of the flow, that 
$$
\int_{\{ABCD\}} \nabla\cdot u\,dx=\int_{\partial
  \{ABCD\}}u\cdot\nv\,d\sigma=0,  
$$
for each (also curvilinear or infinitesimal) ``rectangle'' $\{ABCD\}$ 
touching the boundary $\Gamma$ as in the figure. Since the behaviour
of the flow is not known, in general we have\footnote{Note that the
  line integral over the segment $\{AB\}$ vanishes, since on the
  ``true boundary'' $\Gamma$ both the Navier or no-slip conditions
  impose that $u\cdot\nv=0.$} 
\begin{equation*}
\int_{\{CD\}}u\cdot\nv\,ds
=-\left[\int_{\{BC\}}u\cdot\nv\,ds+
\int_{\{DA\}}u\cdot\nv\,ds\right]\not=0.   
\end{equation*}
This may justify the introduction of a non vanishing normal flux,
also with very low regularity properties, namely the same shared by
the trace of a turbulent flow in the boundary layer region.
\subsection{Main results}
In this section we briefly enunciate the results we will prove,
together with their precise and rigorous statement.

In the sequel $\Omega$ will a bounded, connected, open set in $\R^3,$  
locally situated on one side of its boundary $\Gamma,$ a manifold
of (at least) class $C^{1,1}$ (Lipschitz-continuous first
derivatives). The existence 
of  the unit outward normal $\nv$ derives by results proved in
Ne\v{c}as~\cite{MR37:3168}.

The first result we will prove is an existence and uniqueness theorem
for weak solutions of the Navier-Stokes system~\eqref{maina}, with
boundary conditions~\eqref{normal_component}.
\begin{theorem}\label{mainex}
Let be given $g\in H^{\frac12+\ep}(0,T;H^{\frac12}(\Gamma))$,
satisfying the compatibility condition~\eqref{compa}, $f\in
L^2((0,T)\times\Omega)$, and assume that $u_0\in L^2(\Omega)$, with
$\Div u_0=0$. Then, there exists a unique weak solution
$$
u\in L^\infty(0,T;L^2(\Omega))\cap L^2(0,T;H^1(\Omega)),
$$
of problem~\eqref{maina}.
\end{theorem}

Next, we want to study the behaviour of the solution to
problem~\eqref{maina} as the small parameter $\delta$ converges to
zero. In view of the considerations of the previous section, one can
expect that, as the boundary  layer becomes thinner and thinner, the
solutions will look more and more like the classical solutions
corresponding to the \emph{no-slip} boundary condition. Indeed, this
is the case, as it is shown by the following theorem. \par 
Let $\ud$ be the solution of~\eqref{maina} (we emphasize the
dependence on $\delta$ in this framework) and let $v$ be the solution 
to the Navier-Stokes equations with the same initial value and no-slip
boundary conditions (to be more precise, the vector field $v$ is the
solution to system~\eqref{time_evolution_nse}). As we shall see in
Section~\ref{seclim}, 
$$
\ud=v+\mathcal{O}(\delta^{\frac13}),
$$
so that the ``no-slip solution'' represents the average behaviour,
once one neglects the effect at the boundary. The term $\ud-v$ can be 
seen as the ``fluctuation term'', which takes into account the
non-trivial dynamics at the boundary.
\begin{theorem}\label{mainlim}
Assume $u_0\in H^1(\Omega)$, with $\Div u_0=0$,
$g\in H^{\frac12+\ep}(0,T;H^{\frac12}(\Gamma))$ satisfying the
compatibility condition~\eqref{compa}, and $f\in
L^2((0,T)\times\Omega)$. Then
$$
\sup_{0\le t\le
  T}\|\ud-v\|^2+\int_0^T\Big(\|\D(\ud-v)\|^2+
\frac1\delta\|(\ud-v)\cdot\tv\|^2_\Gamma\Big)\,dt 
=\mathcal{O}(\delta^{\frac23}).
$$
In particular, $\ud$ converges to $v$ in $L^\infty(0,T;L^2(\Omega))$
and $L^2(0,T;H^1(\Omega))$.
\end{theorem}
\section{A result of existence and uniqueness of weak
  solutions}\label{proofex} 
In this section we prove Theorem~\ref{mainex}. For the sake of
simplicity, we consider the normal datum as 
$$
u\cdot\nv=g(x,t),
$$
namely, dropping the dependence on $\delta$, since it is not relevant
in view of the existence and uniqueness result we are going to
show. In the last section we shall see how to deal with a right-hand
side that scales by a power of $\delta$. 

Let us consider the evolution problem
\begin{equation}\label{nse_nonsmooth}
\begin{cases}
{\partial_t}{u}-\nabla\cdot\D(u)+(u\cdot\nabla)\,u+\nabla p=
f&\text{in }\Omega\times(0,T)\\ 
\nabla\cdot u=0&\text{in }\Omega\times(0,T)\\
u\cdot\nv=g(x,t)&\text{on }\Gamma\times(0,T)\\
\delta\,\nv\cdot\D(u)\cdot\tv+u\cdot\tv=0&\text{on }\Gamma\times(0,T)\\  
u(x,0)=u_0(x)&\text{in }\Omega,
\end{cases}
\end{equation}
with $f\in L^2((0,T)\times\Omega)$ and $\Div f=0$, just to avoid
technicalities, and with $g$ not very smooth, say 
\begin{equation}
  \label{g}
g\in H^{1/2+\epsilon}(0,T;H^{1/2}(\Gamma))
\end{equation}
satisfying the compatibility condition~\eqref{compa}.
\subsection{Function spaces}
We use the classical Lebesgue spaces and in particular, we will work
exclusively in the Hilbert framework, so we will use the space
$L^2.$ For simplicity we do not distinguish between space of
scalar, vector, or either tensor valued functions, and the symbol
$\|\,.\,\|$ will denote the norm in $L^2(\Omega).$ The norm in
$L^2(\Gamma)$ will be denoted by $\|\,.\,\|_{\Gamma}.$ \par
In the sequel we will use the customary Sobolev spaces for which we
refer to Adams~\cite{Ada75} and we will use the notion of ``trace''
over the boundary $\Gamma$ of $\Omega.$ Mainly we will use the space
$H^1(\Omega)$ with norm denoted by $\|\,.\,\|_{H^1}$ and its trace
space $H^{1/2}(\Gamma),$ with norm $\|\,.\,\|_{1/2,\Gamma}.$\par 
In addition, we define the spaces
$$
H=\{\,u\in L^2(\Omega)\,|\,\Div u=0,\ u\cdot\nv=0\text{ on
}\Gamma\,\}, 
$$
and
$$
V=\{\,u\in H^1(\Omega)\,|\,\Div u=0,\ u\cdot\nv=0\text{ on
}\Gamma\,\}, 
$$
and we endow $V$ with the norm 
$\|u\|_V=\|\nabla u\|$. Moreover, we define the space of tangential 
vector fields as 
$$
H^1_\tau=\{\,u\in H^1(\Omega)\,|\,u\cdot\nv=0\text{ on }\Gamma\,\}.
$$
\subsubsection{Fractional derivative} In order to define properly the
spaces we will use, we need also to define fractional derivatives. The
fractional derivative may be defined through singular integrals  
$$
D^\alpha_t U(x,t)=\frac{d}{dt}\int_0^t\frac{U(s,x)}{(t-s)^\alpha}\,ds
\quad\text{for }0\leq\alpha<1,
$$
but for our purposes it is better to deal with a definition
\textit{via} the Fourier transform. 

Given $\phi(x,t)$, defined for $t\in[0,T]$,
with values in the Hilbert space $\X$ and integrable (in the Bochner
sense), we define  
\begin{equation}\label{ftilde}
\widetilde{\phi}(t,x)=\begin{cases}
\phi(t,x)&\qquad\text{for }t\in[0,T]\\
0&\qquad\text{elsewhere,}
\end{cases}
\end{equation}
and its Fourier transform (with respect to the time variable) is 
$$
\widehat{\phi}(x,\xi)=\frac{1}{\sqrt{2\pi}}
\int_{\R}\widetilde{\phi}(x,t)\,\e^{-\im t\xi}\,dt,  
$$
so that we can define the fractional Sobolev spaces of functions
having $\alpha-$order derivative in $L^2:$
$$
H^\alpha(\R;\X):=\Big\{f\in L^2(\R;\X):
\ \int_\R{|\xi|^{2\alpha}}\|\widehat{f}(\xi)\|_\X^2
\,d\xi<\infty\Big\}. 
$$
\subsection{The linear stationary problem}
The first step to solve~\eqref{nse_nonsmooth} is to  consider the
linear stationary problem 
\begin{equation}\label{stationary_nonsmooth}
\begin{cases}
-\nabla\cdot\D(G)+\nabla\Pi=0&\text{in }\Omega\times(0,T)\\
\nabla\cdot G=0&\text{in }\Omega\times(0,T)\\
G\cdot\nv=g(x,t)&\text{on }\Gamma\times(0,T)\\
\delta\,\nv\cdot\D(G)\cdot\tv+G\cdot\tv=0&\text{on }\Gamma\times(0,T),
\end{cases}
\end{equation}
where the time variable is now just a parameter.
\begin{theorem}\label{thstationary}
Let be given $g\in H^{1/2+\epsilon}(0,T;H^{1/2}(\Gamma))$, satisfying
the compatibility condition~\eqref{compa}. Then, there exists a unique
$G$ solution of~\eqref{stationary_nonsmooth} such that
\begin{equation}\label{G}
G(x,t)\in H^{1/2+\epsilon}(0,T;H^1(\Omega)).
\end{equation}
Moreover, there is a constant $C_0$, depending only on $\Omega$, such
that 
\begin{equation}\label{stimaG}
\|\nabla G\|+\|\Pi\|\le
C_0(1+\delta^{-\frac12})\|g\|_{H^{\frac12}(\Gamma)}. 
\end{equation}
\end{theorem}
\begin{proof}
See Solonnikov and {\v{S}}{\v{c}}adilov~\cite{MR51:1164} and
Beir\~ao da Veiga~\cite{bdv2004}. In fact, for each $t$ it is possible 
to solve a linear stationary Stokes problem (with the appropriate
boundary conditions) that has a unique solution belonging to
$H^1(\Omega)$. The regularity in the time variable is inherited by the
function $G$. 

We give a formal (but completely justified) argument for the
estimate~\eqref{stimaG}, following the approach to the existence in
$L^2$-spaces introduced in Beir\~ao da Veiga~\cite{bdv2004} to find
appropriate estimates on $G$. We consider the bilinear form
$$
B(u,\phi)=\int_\Omega\D(u)\D(\phi)\,dx
$$ 
and the functions $(G,\Pi),$ that solve~\eqref{stationary_nonsmooth},
must satisfy 
$$
B(G,\phi)-\int_\Omega
\Pi\,\nabla\cdot\phi\,dx+\frac{1}{\delta}\int_\Gamma 
G\cdot\phi\,d\sigma=0,
 \qquad\forall\,\phi\in H^1_\tau.
$$
In order to deal with the inhomogeneous problem, we introduce a vector 
field $G_1$ such that
$$
\begin{cases}
\nabla \cdot G_1=0&\text{ in }\Omega,\\ 
G_1\cdot\nv=g&\text{ on }\Gamma,\\
\|G_1\|_{H^1}\leq C\|g\|_{H^{\frac12}(\Gamma)}.
\end{cases}
$$
The construction of such a vector field is rather standard and can be
found for instance in Galdi~\cite{Gal94}.

By defining $G=G_1+G_2,$ the function $G_2$ must satisfy
$$
B(G_2,\phi)-\int_\Omega\Pi\,\nabla\cdot\phi\,dx+
\frac{1}{\delta}\int_\Gamma  G_2\cdot\phi\,d\sigma
=-B(G_1,\phi)-\frac{1}{\delta}\int_\Gamma G_1\cdot\phi\,d\sigma,
$$
for each $\phi$ tangential to the boundary. If $\phi=G_2$, we get
(since $G_1\cdot\nv=0$ and $\nabla\cdot G_2=0$)
$$
\|\nabla G_2\|^2+\frac{1}{\delta}\|G_2\|_\Gamma^2
\leq
\|\nabla G_2\|\,\|\nabla
G_1\|+\frac{1}{\delta}\|G_1\|_\Gamma\|G_2\|_\Gamma,  
$$
and consequently 
$$
\frac12\|\nabla G_2\|^2+\frac1{2\delta}\|G_2\|_\Gamma^2
\le\frac12\|\nabla G_1\|^2+\frac1{2\delta}\|G_1\|_\Gamma^2. 
$$
This finally implies that
$$
\|\nabla G_2\|^2\leq
C\left(1+\frac1\delta\right)\|g\|_{H^{\frac12}(\Gamma)}^2, 
$$
where the constant $C$ depends on $\Omega,$ but it is independent of
$\delta$. The estimate on the pressure can be obtained by
approximation, by studying a slightly different equation
(see~\cite{bdv2004} for details).  
\end{proof}

Indeed, notice that, for our aims, it is enough to have
\begin{equation}\label{assG}
G(x,t)\in H^{1/2+\epsilon}(0,T;L^2(\Omega))\cap
L^{2}(0,T;H^1(\Omega)).  
\end{equation}
but currently we do not know the minimal assumption on $g$
in order to have the above regularity. Regarding the usual no-slip 
boundary conditions see the result proved by Fursikov,
Gunzburger, and Hou~\cite{MR2002j:35244}.
\subsection{The linear evolution problem}
The next step for the analysis of the nonlinear evolution
problem~\eqref{nse_nonsmooth} is the following linear evolution
problem 
\begin{equation}\label{se_nonsmooth}
\begin{cases}
{\partial_t}{z}-\nabla\cdot\D(z)+\nabla q=0&\text{in
}\Omega\times(0,T)\\ 
\nabla\cdot z=0&\text{in }\Omega\times(0,T)\\
z\cdot\nv=g&\text{on }\Gamma\times(0,T)\\
\delta\,\nv\cdot\D(z)\cdot\tv+z\cdot\tv=0&\text{on }\Gamma\times(0,T)\\
z(x,0)=G(x,0)&\text{in }\Omega.
\end{cases}
\end{equation}
We shall treat the non-linear problem as a perturbation of such a
linear system. Let us introduce the new unknowns
$$
Z(x,t)=z(x,t)-G(x,t)\quad\text{and}
\quad Q(x,t)=q(x,t)-\Pi(x,t)
$$
so that we are reduced to a homogeneous problem for the new unknowns
$(Z,Q)$: 
\begin{equation}\label{Vse_nonsmooth}
\begin{cases}
{\partial_t}{Z}-\nabla\cdot\D(Z)+\nabla Q=-\partial_t G&\text{in
  }\Omega\times(0,T)\\
 \nabla\cdot Z=0&\text{in }\Omega\times(0,T)\\
Z\cdot\nv=0&\text{on }\Gamma\times(0,T)\\
\delta\,\nv\cdot\D(Z)\cdot\tv+Z\cdot\tv=0&\text{on
  }\Gamma\times(0,T)\\ 
Z(x,0)=0&\text{in }\Omega,
\end{cases}
\end{equation}
The above problem is not completely standard, since the right-hand
side does not satisfy the usual properties. For instance,
one can note that $\partial_t G$ does not belong to the domain of the
Stokes operator, since $\partial_t G\cdot\nv\not=0$. This is the main
difficulty: the low regularity of this term can be treated in a more
standard way, while the above fact is responsible for a different
approach.
\begin{theorem}\label{thlinear}
Assume that $(G,\Pi)$ is a solution to
system~\eqref{stationary_nonsmooth}, with $G$ satisfying the
regularity property~\eqref{assG}. Then, there exists a unique solution
$(z,q)$ (where $q$ is unique up to an additive constant not depending
on the space variable) to system~\eqref{se_nonsmooth} such that
\begin{equation}\label{regz}
z\in L^\infty(0,T;L^2(\Omega))\cap L^2(0,T;H^1(\Omega))\cap
H^{\frac12-\ep}(0,T;L^2(\Omega)). 
\end{equation}
Moreover,
\begin{multline}\label{zestim}
 \sup_{0\le t\le T}\|z(t)\|^2
+\int_0^T\Big(\|\D(z)(t)\|^2+\frac1\delta\|z(t)\cdot\tv\|^2_\Gamma\Big)\,dt 
+\|z\|^2_{H^{\frac12-\ep}(0,T;L^2(\Omega))}\le\\
\le C\left(\|G\|^2_{H^{\frac12+\ep}(0,T;L^2(\Omega))}
         +\int_0^T\|\D(G)\|^2\right).
\end{multline}
\end{theorem}
\begin{proof}
By virtue of Theorem~\ref{thstationary}, it is enough to prove the
same claim of this theorem on the solution $(Z,Q)$ of
problem~\eqref{Vse_nonsmooth}. Since we just know that $\partial_t G\in
H^{-\frac12+\ep}(0,T;L^2(\Omega)),$ we introduce a sequence $G^N\in
H^1(\R;L^2(\Omega))$ of approximate functions such that
\begin{itemize}
\item[\textit{(a)}] $G^N|_{[0,T]}\longrightarrow G$ in
  $H^{\frac12+\ep}(0,T;L^2(\Omega)),$ as $\to\infty$,
\item[\textit{(b)}] $\|\partial_t
  G^N\|_{L^2(0,T;L^2(\Omega))}={N}$.
\end{itemize}
The way to do this extension is rather standard: first we can define
$\overline{G}:\,\R\to L^2(\Omega)$ with an extension by reflection. 
Then, we consider a sequence $\rho_N$ of mollifiers 
and the function $G_N$ will be the restriction on $[0,T]$ of the
function $\rho_N*\overline{G}$. 

The proof is based on the Faedo-Galerkin procedure. By
Clopeau~\textit{et al.}~\cite{MR99g:35102}, (but also the recent
abstract results in~\cite{bdv2004}) we know that there exists a basis 
$(\phi_n)_{n\in\N}$ of functions in $H^3(\Omega)$ of the space $V$
(and also of $H$), such that 
$$
\delta\,\nv\cdot\D(\phi_n)\cdot\tv+\phi_n\cdot\tv=0.
$$
Now, let
$Z_n^N(t,x)=\sum_{k=1}^n\zeta_{n,k}^N(t)\phi_k(x)$ be the solution of
the following (finite dimensional) linear system of Ordinary
Differential Equations (ODE):
$$
\begin{cases}
\displaystyle{\frac{d}{dt}
  \langle Z_n^N,\phi_k\rangle
 +\langle\D(Z_n^N),\D(\phi_k)\rangle
 +\frac{1}{\delta}\int_\Gamma(Z_n^N\cdot\tv)\,(\phi_k\cdot\tv)\,d\sigma
=-\frac{d}{dt}\langle G^N,\phi_k\rangle}\\ 
(Z_n^N(x,0),\phi_k)=0,&
\end{cases}
$$
for $t\in(0,T)$ and $k=1,\dots,n$,
where
$\langle\cdot,\cdot\rangle$ denotes the usual
$L^2$-scalar product. Notice that the divergence-free
constraint and the boundary conditions on $Z_n^N$ are automatically
verified. By using a standard argument it is easy to prove that such a
system of ODE has a unique solution $Z_n^N\in
L^\infty(0,T;L^2(\Omega))\cap L^2(0,T;H^1(\Omega));$ indeed, by
multiplying each equation by the corresponding term $\zeta^N_{n,k}$,
summing over $k$ and integrating by parts over $\Omega$, one
easily obtains the following estimate  
$$
\sup_{0\le t\le T}\|Z^N_n(t)\|^2+\int_0^T\Big(\|\D(Z^N_n)\|^2+
\frac1\delta\int_0^T\|Z^N_n\|^2_\Gamma\Big)\,dt 
\le C\|G^N\|^2_{H^1(0,T;L^2(\Omega))},
$$
for a constant $C$, depending only on $\Omega$. Unfortunately, such an
estimate, beside being uniform in $n$, is not uniform in $N$ due to
property \textit{(b)} of the approximate sequence
$(G^N)_{N\in\N}.$ Hence, we need other \emph{a-priori} estimates on
the solutions $Z_n^N$ of the finite dimensional problem. We continue
working on the $Z_n^N$, since such functions are smooth enough for the
computation that will be performed.

We multiply again the equations by the terms $\zeta^N_{n,k}$, sum
over $k$ and integrate by parts, but we 
estimate the right-hand side in the following way:
\begin{equation}\label{fracest}
\begin{split}
\sup_{0\le t\le T}\|Z^N_n(t)\|^2
&+\int_0^T\left(\|\D(Z^N_n)\|^2
 +\frac1\delta
\|Z^N_n\|^2_\Gamma\right)\,dt
\le\Big|\int_0^T\int_\Omega\partial_tG^n\cdot Z_n^N\Big|\,dxdt\\
&\le \|\partial_tG^N\|_{H^{-\frac12+\ep}(0,T;L^2(\Omega))}
\|Z_n^N\|_{H^{\frac12-\ep}(0,T;L^2(\Omega))}\\ 
&\le\|G^N\|_{H^{\frac12+\ep}(0,T;L^2(\Omega))}
\|Z_n^N\|_{H^{\frac12-\ep}(0,T;L^2(\Omega))}, 
\end{split}
\end{equation}
so that we only need to show a uniform estimate (with respect to both
$n$ and $N$) of $Z_n^N$ in the space
$H^{\frac12-\ep}(0,T;L^2(\Omega)).$ We shall use the Fourier transform 
characterization of the norm of fractional Sobolev spaces (see
Adams~\cite{Ada75}) to get such estimate. Each $\tZ$ (such functions
have been defined in~\eqref{ftilde}) is a solution of the following
equation
\begin{multline*}
  \frac{d}{dt}\int_\Omega\tZ\cdot\phi_k
 +\int_\Omega\D(\tZ)\cdot\D(\phi_k)
 +\frac1\delta\int_\Gamma\tZ\cdot\phi_k=\\
=-\frac{d}{dt}\int_\Omega\widetilde{G}^N\cdot\phi_k
 +\delta(t)\int_\Omega G^N(0)\cdot\phi_k
 -\delta(t-T)\int_\Omega(Z_n^N(T)+G^N(T))\cdot\phi_k,
\end{multline*}
for each $k=1,\ldots,n,$ in the sense of distributions with respect to
the time variable. Here $\delta(\cdot)$ is the usual Dirac's
\emph{delta function}. In the frequency Fourier variable $\xi,$ the
above equation reads: 
\begin{multline*}
 -\im\xi\int_\Omega\hZ\cdot\phi_k
 +\int_\Omega\D(\hZ)\cdot\D(\phi_k)
 +\frac1\delta\int_\Gamma\hZ\cdot\phi_k=\\
=\im\xi\int_\Omega\widehat{G}^N\cdot\phi_k
 +\int_\Omega G^N(0)\cdot\phi_k
 -\e^{-\im\xi T}\int_\Omega(Z_n^N(T)+G^N(T))\cdot\phi_k,
\end{multline*}
see for instance Lions~\cite{MR41:4326}, where this tool is used to
prove estimates on the fractional derivative of the solution. Note
that in that reference, and in all involving fractional derivatives
for the Navier-Stokes equations, the starting point is the existence
of a weak solution and on it it is possible to prove additional
estimates. In our case the existence of weak solution derives from the
fractional derivative estimates and it seems not possible to prove the
usual existence results. 

Consequently, we get
\begin{multline*}
 -\im\xi\|\hZ(\xi)\|^2
 +\|\D(\hZ)(\xi)\|^2
 +\frac1\delta\|\hZ(\xi)\|_\Gamma^2=\\
=\im\xi\int_\Omega\widehat{G}^N\cdot\overline{\hZ}
 +\int_\Omega G^N(0)\cdot\overline{\hZ}
 -\e^{-\im\xi T}\int_\Omega(Z_n^N(T)+G^N(T))\cdot\overline{\hZ}.
\end{multline*}
Take the imaginary part and multiply both sides of the previous formula by
$|\xi|^{2\alpha-1}$, with $\alpha<\frac12$ so that, by using Young's
inequality, one gets
$$
|\xi|^{2\alpha}\|\hZ(\xi)\|^2
\le C|\xi|^{2\alpha}\|\widehat{G}^N\|^2
   +C|\xi|^{2\alpha-2}(\|G^N(T)\|+\|Z_n^N(T)\|+\|G^N(0)\|)^2.
$$
In order to estimate the integral
$\int_\R|\xi|^{2\alpha}\|\hZ(\xi)\|^2\,d\xi$, 
we split it in two parts; by the above estimate,
\begin{multline*}
\int_{|\xi|>1}|\xi|^{2\alpha}\|\hZ(\xi)\|^2\le\\
\le C\int_\R|\xi|^{2\alpha}\|\widehat{G}^N\|^2
   +C(\|G^N(T)\|+\|Z_n^N(T)\|+\|G^N(0)\|)^2\int_{|\xi|>1}|\xi|^{2\alpha-2};
\end{multline*}
the first term on the right-hand side is controlled by
$C\|G^N\|^2_{H^{\frac12+\ep}(0,T;L^2(\Omega))}$, while
$$
\|Z_n^N(T)\|^2\le
C\|G^N\|_{H^{\frac12+\ep}(0,T;L^2(\Omega))}
\|Z_n^N\|_{H^{\frac12-\ep}(0,T;L^2(\Omega))}, 
$$
by virtue of~\eqref{fracest}; $\|G^N(0)\|$ is bounded by $\|G(0)\|$,
and finally, by using the Morrey inequality that implies
$H^{1/2+\ep}(0,T)\subset 
C([0,T])$ (see Adams~\cite{Ada75}) we get 
$$
\|G^N(T)\|\le\|G^N\|_{H^{\frac12+\ep}(0,T;L^2(\Omega))}.
$$
The second part is estimated as follows, by using Parseval's theorem,
Poincar\'e inequality and estimate~\eqref{fracest},
\begin{align*}
\int_{|\xi|\le1}|\xi|^{2\alpha}\|\hZ\|^2\,d\xi
&\le\int_\R\|\hZ\|^2\,d\xi
=\int_0^T\|Z_n^N(t)\|^2\,dt\\
&\le C\int_0^T\|\D(Z^N_n)\|^2dt\\
&\le
C\|G^N\|_{H^{\frac12+\ep}(0,T;L^2(\Omega))}\|Z_n^N\|_{H^{\frac12-\ep}(0,T;L^2(\Omega))}.  
\end{align*}
In conclusion, by collecting all of the above estimates, we finally
get that, for each $\ep\in(0,\frac12)$, there exists a constant $C$,
depending only on $\Omega$ and $\ep$, such that
\begin{equation}\label{fracest2}
\|Z_n^N\|_{H^{\frac12-\ep}(0,T;L^2(\Omega))}
\le C\|G^N\|_{H^{\frac12+\ep}(0,T;L^2(\Omega))},
\end{equation}
that, together with~\eqref{fracest}, says that $Z_n^N$ is bounded,
uniformly in $n$ 
and $N$, in the spaces $H^{\frac12-\ep}(0,T;L^2(\Omega))$,
$L^\infty(0,T;L^2(\Omega))$ 
and $L^2(0,T;H^1(\Omega))$. 

Hence, it is possible to extract a (diagonal) subsequence
converging weakly in $L^2(0,T;H^1(\Omega))$, \mbox{weakly-$*$} in
$L^\infty(0,T;L^2(\Omega))$ 
and strongly in $L^2((0,T)\times\Omega)$ to the unique solution $Z$ of problem
\eqref{Vse_nonsmooth}. Indeed, $Z\in
H^{\frac12-\ep}(0,T;L^2(\Omega))$, that is, the 
topological dual space of $H^{-\frac12+\ep}(0,T;L^2(\Omega))$, the
space to which 
$\partial_t G$ belongs to. Furthermore, by passing to the limit and by
using the 
semi-continuity of the norms, it follows that $Z$ satisfies the claim
stated at the beginning of this proof.
\end{proof}
\begin{remark}
The assumption $g\in H^{\frac12+\ep}(0,T;H^{\frac12}(\Gamma))$, which
in turn gives 
$G\in H^{\frac12+\ep}(0,T;H^1(\Omega))$, seems to be rather technical, for the
presence of $\epsilon$. If $g\in
H^{\frac12}(0,+\infty;H^{\frac12}(\Gamma))$, it 
follows that $G\in H^{\frac12}(0,+\infty;H^1(\Omega))$ and
Theorem~\ref{thlinear}  holds accordingly.

Indeed, the main point is estimate~\eqref{fracest}, in which the
right-hand side becomes $\|G^N\|_{H^{\frac12}}\|Z_n^N\|_{H^{\frac12}}$ and,
following the lines of the proof presented above, the estimate on the
Fourier transform gives that $Z^N_n$ is bounded, uniformly in $n$ and
$N$, in $H^{\frac12}(0,+\infty;L^2(\Omega))$. Notice that in the
\emph{critical} case $H^{\frac12}$, we work on the whole time interval
$[0,+\infty)$, to avoid the boundary terms $\|G^N(0)\|$ and
$\|G^N(T)\|$, that cannot be estimated by using the Morrey inequality.

We also note that this small relaxation on the assumptions on $G$
requires to add the
hypothesis $G\in L^\infty(0,T;L^2(\Omega)).$ Otherwise the function $z$
will not belong itself to $L^\infty(0,T;L^2(\Omega))$ and this fact
is crucial to prove the corresponding bound for weak solutions to the
full nonlinear Navier-Stokes problem. 
\end{remark}
\subsection{The non-linear problem}
In this section we finally prove Theorem~\ref{mainex}. Again, we make use
of an auxiliary problem; namely, we introduce the new variables
$$
U=u-z\qquad\text{and}\qquad P=p-q,
$$
where $(z,q)$ is the solution to the linear evolution
problem~\eqref{se_nonsmooth}, and the pair $(U,P)$ solves the
following problem 
\begin{equation}\label{nse_homo}
\begin{cases}
\partial_t U-\nabla\cdot\D(U)+[(U+z)\cdot\nabla](U+z)+\nabla
P=f&\text{in }\Omega\times(0,T)\\ 
\nabla\cdot U=0&\text{in }\Omega\times(0,T)\\
U\cdot\nv=0&\text{on }\Gamma\times(0,T)\\
\delta\,\nv\cdot\D(U)\cdot\tv+U\cdot\tv=0&\text{on
}\Gamma\times(0,T)\\ 
U(x,0)=u_0(x)-G(x,0)&\text{in }\Omega.
\end{cases}
\end{equation}
By virtue of Theorem~\ref{thlinear}, the existence Theorem~\ref{mainex}
for the non-linear problem is a straightforward consequence of the
following proposition.
\begin{proposition}
Assume that $(G,\Pi)$ is solution to system~\eqref{stationary_nonsmooth},
with $G\in H^{\frac12+\ep}(0,T;L^2(\Omega))\cap L^2(0,T;H^1(\Omega))$,
then there is a unique
\begin{equation}\label{regU}
U\in L^\infty(0,T;L^2(\Omega))\cap L^2(0,T;H^1(\Omega))
\end{equation}
solution to problem~\eqref{nse_homo}. Moreover, the following estimate 
holds true:
\begin{multline}\label{Uestim}
\sup_{0\le s\le
  t}\|U(s)\|^2+\int_0^t\Big(\|\D(U)\|^2+
\frac1\delta\|U\|_\Gamma^2\Big)\,ds\le\\  
\le \|u_0-G(\cdot,0)\|^2\e^{A(t)}+C\int_0^t(\|f\|^2+\|\nabla
  z(s)\|^2\|z(s)\|)
\e^{A(t)-A(s)}\,ds,
\end{multline}
where 
$$
A(t)=Ct+C\left(1+\|z\|^2_{L^\infty(0,T;L^2(\Omega))}\right)
\int_0^t\|\nabla z\|^2\,ds,
$$
and $C$ is a constant depending only on $\Omega$.
\end{proposition}
\begin{proof}
The proof is rather standard and proceeds \textit{via} a
Faedo-Galerkin approximation, as in the proof of
Theorem~\ref{thlinear}. We only show an \emph{a-priori} estimate,
whose computations are formal, but are completely meaningful at the
level of the Faedo-Galerkin approximate functions. Multiply
\eqref{nse_homo} by $U$ and integrate by parts, to get 
$$
\frac12\frac{d}{dt}\|U\|^2+\|\D(U)\|^2+\frac1\delta\|U\|^2_\Gamma
=\int_\Omega U\cdot[(U+z)\cdot\nabla](U+z)+\int_\Omega f\cdot U.
$$
The estimate of the integral involving $f$ is straightforward, since
it is bounded by $\|f\|^2+\|U\|^2.$ We estimate the non-linear term in
the right-hand side by using the Gagliardo-Nirenberg
inequality\footnote{Note that such a inequality 
is a little bit more general than the so-called
\emph{Lady\v{z}henskaya inequality}, since the functions are not
vanishing on the boundary of $\Omega$.}, 
\begin{equation}\label{g-n}
\|u\|_{L^4}\leq C\|u\|^{1/2}\|\nabla u\|^{1/2}\qquad\forall\,u\in
H^1(\Omega). 
\end{equation}
We first observe that, since $\Div U=0$ and $U\cdot\nv=0$,
$$
\int_\Omega U\cdot(U\cdot\nabla)U=0
\qquad\text{and}\qquad
\int_\Omega U\cdot(U\cdot\nabla)z=-\int_\Omega z\cdot(U\cdot\nabla)U,
$$
so that, using repeatedly the Gagliardo-Nirenberg inequality above
given, H\"older's inequality and Young's inequality, we get
\begin{align*}
\left|\int_\Omega U\cdot[(U+z)\cdot\nabla](U+z)\right|
&\le 2\|U\|_{L^4}\|z\|_{L^4}\|\nabla U\|
    +\|U\|_{L^4}\|z\|_{L^4}\|\nabla z\|\\
&\le \frac12\|\nabla U\|^2+C\|z\|^4_{L^4}\|U\|^2
    +C\|z\|^{\frac43}_{L^4}\|\nabla z\|^{\frac43}\|U\|^{\frac23}\\
&\le \frac12\|\nabla U\|^2+C\|z\|^2\|\nabla z\|^2\|U\|^2
    +C\|z\|^{\frac23}\|\nabla z\|^2\|U\|^{\frac23}\\
&\le \frac12\|\nabla U\|^2+C(1+\|z\|^2)\|\nabla z\|^2\|U\|^2
    +C\|\nabla z\|^2\|z\|.
\end{align*}
Since, by Theorem~\ref{thlinear}, $z\in L^\infty(0,T;L^2(\Omega))\cap
L^2(0,T;H^1(\Omega))$, both terms $(1+\|z\|^2)\|\nabla z\|^2$ and
$\|\nabla z\|^2\|z\|$ are integrable in time and, by Gronwall's lemma,
we can deduce that $U$ is bounded in $L^\infty(0,T;L^2(\Omega))$ and
$L^2(0,T;H^1(\Omega))$. Moreover, formula~\eqref{Uestim} also
follows.\par 
Finally, uniqueness of the solution follows from similar
arguments. Indeed, if $\widetilde{U}$ is the difference between two
solutions $U_1$ and $U_2,$ one easily gets
$$
\frac12\frac{d}{dt}\|\widetilde{U}\|^2
\le(\|U_2\|_{L^4}^4+\|z\|_{L^4}^4+\|\nabla z\|^2)\|\widetilde{U}\|^2
$$
and, since $\widetilde{U}(0)=0$, from Gronwall's lemma it follows that
$\widetilde{U}\equiv0$.
\end{proof}
\begin{remark}\label{2D-3D}
  In the proof of the result of this section we used in a fundamental
  way estimate~\eqref{g-n}. In the three-dimensional case this
  inequality is not anymore true. Instead, it holds
\begin{equation*}
\|u\|_{L^4}\leq C\|u\|^{1/4}\|\nabla u\|^{3/4}\qquad\forall\,u\in
H^1(\Omega), 
\end{equation*}
and the latter can be used to prove just local existence of weak
solutions. The global result proved in the 2D case stems in an
essential manner on the stronger estimate and this is the critical
difference between the two cases.
\end{remark}
\section{The proof of Theorem~\ref{mainlim}}\label{seclim}
In this section we study the limit of the solution of the vorticity
seeding model as $\delta\to0.$ We prove the convergence of
the solutions of~\eqref{time_evolution_delta} to the solutions of the
non-stationary Navier-Stokes system with the no-slip boundary
condition~\eqref{time_evolution_nse}. This supports the idea of using
a non-standard boundary condition on the new boundary
$\Gamma_1\subset\Omega,$ such that  the region between $\Gamma$ and
$\Gamma_1$ is very narrow. When the width of this region (presumably
the boundary layer) shrinks to zero, then the solution of the
Navier-Stokes equations with the usual no-slip boundary condition is
recovered.   
\subsection{Comparison of solutions with different boundary data}
Let $g\in H^{\frac12+\ep}(0,T;H^{\frac12}(\Gamma))$, satisfying the
compatibility condition~\eqref{compa}, $f\in L^2((0,T)\times\Omega)$
and $u_0\in H^1(\Omega)$, with $\Div u_0=0$. Without loss of generality,
we can assume that $u_0\equiv0$. Denote by $\Gd$ the
solution of the linear stationary problem~\eqref{stationary_nonsmooth},
with boundary condition $\Gd\cdot\nv=\delta\, g$. Consider the solution
$(\ud,p_\delta)$ to the system
\begin{equation}\label{time_evolution_delta}
\begin{cases}
\partial_t\ud-\nabla\cdot\D(\ud)+\ud\cdot\nabla\ud+\nabla p_\delta=f
  &\text{in }\Omega\times(0,T)\\ 
\nabla\cdot\ud=0&\text{in }\Omega\times(0,T)\\
\ud\cdot\nv=\delta\, g(x,t)&\text{on }\Gamma\times(0,T)\\
\delta\,\nv\cdot\D(\ud)\cdot\tv+\ud\cdot\tv=0
  &\text{on }\Gamma\times(0,T)\\ 
\ud(x,0)=u_0(x)+G(x,0)&\text{in }\Omega,
\end{cases}
\end{equation}
We want to show the convergence of the vector valued function $\ud$ to 
the solution $v$ of the Navier-Stokes equations with zero Dirichlet
data:  
\begin{equation}\label{time_evolution_nse}
\begin{cases}
\partial_t v-\nabla\cdot\D(v)+v\cdot\nabla v+\nabla \pi=f&\text{in
  }\Omega\times(0,T)\\  
\nabla\cdot v=0&\text{in }\Omega\times(0,T)\\
v=0&\text{on }\Gamma\times(0,T)\\
v(x,0)=u_0&\text{in }\Omega.
\end{cases}
\end{equation}
To this end, we also introduce the solution $\zd$ to the linear
evolution problem~\eqref{se_nonsmooth}, with boundary condition
$\zd\cdot\nv=\delta\, g$, and we set
$$
\Ud=\ud-\zd\qquad\text{and}\qquad \od=\ud-\zd-v.
$$
The function $\od$ satisfies the following homogeneous system
\begin{equation}\label{weq}
\begin{cases}
\partial_t\od-\nabla\cdot\D(\od)+R(\od,\zd,v,\Ud)+\nabla r=0
  &\text{in }\Omega\times(0,T)\\ 
\nabla\cdot\od=0&\text{in }\Omega\times(0,T)\\
\od\cdot\nv=0&\text{on }\Gamma\times(0,T)\\
\delta\,\nv\cdot\D(\od)\cdot\tv+\od\cdot\tv=
-\delta\,\nv\cdot\D(v)\cdot\tv 
  &\text{on }\Gamma\times(0,T)\\
\od(x,0)=0&\text{in }\Omega,
\end{cases}
\end{equation}
where
$$
R(\od,\zd,v,\Ud)=(\Ud\cdot\nabla)\od+(\od\cdot\nabla)v
                +(\Ud\cdot\nabla)\zd+(\zd\cdot\nabla)\Ud
                +(\zd\cdot\nabla)\zd.
$$
We multiply the first equation in~\eqref{weq} by $\od$ and integrate
by parts to get 
$$
\frac12\frac{d}{dt}\|\od\|^2+\|\D(\od)\|^2+\frac1\delta\|\od\|_\Gamma^2
=-\int_\Gamma(\nv\cdot\D(v)\cdot\tv)(\od\cdot\tv)d\sigma-
\int_\Omega\od\cdot 
R(\od,\zd,v,\Ud). 
$$
The boundary integral in the right-hand side may be increased in the
following way, 
\begin{equation}\label{boundary}
\left|\int_\Gamma(\nv\cdot\D(v)\cdot\tv)(\od\cdot\tv)\,d\sigma\right|
\le\frac{\delta}{2}\|v\|_{H^2}^2+
\frac{1}{2\delta}\|\od\|_\Gamma^2. 
\end{equation}
We analyse then the second integral. By integration by parts over
$\Omega$ (note that $\Ud\cdot\nv=0$ on $\Gamma$) and by using
H\"older's inequality, the Gagliardo-Nirenberg inequality, and Young's
inequality, we obtain 
\begin{eqnarray*}
\lefteqn{\left|\int_\Omega\od\cdot R(\od,\zd,v,\Ud)\,dx\right|\le}\\
&\le&\frac12\|\D(\od)\|^2
     +C\big[\|\nabla v\|^2+\|\nabla\Ud\|^{\frac23}\|\nabla\zd\|^{\frac43}
     +\|\nabla\Ud\|^{\frac43}\|\nabla\zd\|^{\frac23}+
\|\nabla\zd\|^2\big]\|\od\|^2\\ 
&   &+C\big[\|\nabla\Ud\|^{\frac23}\|\nabla\zd\|^{\frac43}\|\Ud\|
          +\|\nabla\Ud\|^{\frac43}\|\nabla\zd\|^{\frac23}\|\zd\|
          +\|\nabla\zd\|^2\|\zd\|\big].
\end{eqnarray*}
Indeed, for example,
\begin{align*}
\left|\int_\Omega\od\cdot(\zd\cdot\nabla)\Ud\right|
&\le\|\od\|_{L^4}\|\zd\|_{L^4}\|\nabla U\|\\
&\le
C\|\od\|^{\frac12}\|\D(\od)\|^{\frac12}\|\zd\|^{\frac12}
\|\D(\zd)\|^{\frac12}\|\nabla\Ud\|\\ 
&\le\frac18\|\D(\od)\|^2
   +C\|\nabla\Ud\|^{\frac43}\|\nabla\zd\|^{\frac23}\|\od\|^{\frac23}
\|\zd\|^{\frac23}\\
&\le\frac18\|\D(\od)\|^2
   +C\|\nabla\Ud\|^{\frac43}\|\nabla\zd\|^{\frac23}(\|\od\|^2+\|\zd\|).
\end{align*}
For simplicity we set
\begin{align*}
\psi(t)&:=C\big[\|\nabla\Ud\|^{\frac23}\|\nabla\zd\|^{\frac43}\|\Ud\|
          +\|\nabla\Ud\|^{\frac43}\|\nabla\zd\|^{\frac23}\|\zd\|
          +\|\nabla\zd\|^2\|\zd\|\big],\\
\phi(t)&:=\big[\|\nabla v\|^2+\|\nabla\Ud\|^{\frac23}\|\nabla\zd\|^{\frac43}
      +\|\nabla\Ud\|^{\frac43}\|\nabla\zd\|^{\frac23}+
\|\nabla\zd\|^2\big],\\
\chi(t)&:=\|v\|_{H^2}^2,
\end{align*}
so that, by recollecting all together the estimates we have before
obtained, we end up with the following differential inequality 
\begin{equation}\label{limest}
\frac{d}{dt}\|\od\|^2+\|\D(\od)\|^2+\frac1\delta\|\od\|_\Gamma^2
\le\delta\chi(t)+\psi(t)+\phi(t)\|\od\|^2
\end{equation}
\subsubsection{Estimate of $\psi$ and $\phi$ in terms of $\delta$}
First, we note that the functions $\psi$ and $\phi$ belong to
$L^1(0,T)$, by virtue of the regularity properties~\eqref{regz}
and~\eqref{regU} of $\zd$ and $\Ud$, respectively.
The function $\chi(t)$ belongs to $L^1(0,T)$ as well, since
$$
\|v\|_{L^2(0,T;H^2)}\le
C[\|v(x,0)\|_{H^1}+\|f\|_{L^2((0,T)\times\Omega)}], 
$$
see for instance Constantin and Foia{\c{s}}~\cite{CF88}. And also in 
this point is crucial to consider the 2D problem, since such estimate
is available just for small times in the 3D case.

It is important now to sharply check the correct behaviour of the
functions $\phi$ and $\psi(t)$ in terms of $\delta$. Before going
further, we collect in the following lemma the results of
Section~\ref{proofex} that we shall need.
\begin{lemma}
Under the assumptions of Theorem~\ref{mainlim}, the following
estimates hold
\begin{align*}
\textit{(i)}&\quad \sup_{0\le t\le T}\|\zd\|^2+\int_0^T\|\D(\zd)\|^2\,dt
\le C_1\delta\|g\|^2_{H^{\frac12+\ep}(0,T;H^{\frac12}(\Gamma))},\\
\textit{(ii)}&\quad\sup_{0\le t\le T}\|\Ud\|^2+\int_0^T\Big(\|\D(\Ud)\|^2
+\frac1\delta\|\Ud\|^2_\Gamma\Big)\,dt\le C_2,
\end{align*}
for all $0<\delta\le1$, where the constant $C_1$ depends only on
$\Omega$ and the constant $C_2$ depends only on $\Omega$, $T$,
$\|f\|_{L^2}$, and $\|g\|_{H^{\frac12+\ep}}$.
\end{lemma}
\begin{proof}
Property \textit{(i)} is a consequence of~\eqref{zestim}. Indeed,
by~\eqref{stimaG} and Poincar\'e inequality,
$$
\|\Gd\|\le C\|\D(\Gd)\|\le C\sqrt{\delta}\|g\|,
$$
so that $\|\Gd\|_{H^{\frac12+\ep}(0,T;L^2(\Omega))}\le
C\sqrt{\delta}\|g\|_{H^{\frac12+\ep}(0,T;L^2(\Omega))}.$ 

Again, from~\eqref{stimaG},
$$
\int_0^T\|\D(\Gd)\|^2\,dt\le
C\delta\|g\|^2_{H^{\frac12+\ep}(0,T;L^2(\Omega))},
$$
and in conclusion \textit{(i)} holds true. 

As it concerns \textit{(ii)}, we have that, by Poincar\'e inequality,
\eqref{stimaG} and Sobolev embeddings,  
$$
\|\Gd(\cdot,0)\|
\le C\|\D(\Gd)(\cdot,0)\|
\le C\sqrt{\delta}\|g(\cdot,0)\|
\le C\sqrt{\delta}\|g\|_{H^{\frac12+\ep}(0,T;L^2)}.
$$
Moreover, since $\delta\le1$ and from the above estimates it follows 
$$
A(t)\le A(T)\le CT+C(1+\|g\|^4_{H^{\frac12+\ep}(0,T;L^2(\Gamma))})
$$
that is, $A(t)$ is uniformly bounded by a constant independent of
$\delta$, and 
\begin{align*}
\int_0^t(\|f\|^2+\|\nabla\zd\|^2\|\zd\|)\e^{A(t)-A(s)}\,ds
&\le\e^{A(T)}(\int_0^T\|f\|^2+\|\zd\|_{L^\infty(0,T;L^2)}\|\D(\zd)\|^2)\\
&\le(\|f\|^2_{L^2((0,T)\times\Omega)}+C\|g\|^3_{H^{\frac12+\ep}})\e^{A(T)}
\end{align*}
as well, so that, by~\eqref{Uestim}, also \textit{(ii)} follows.
\end{proof}

The first consequence of the above lemma is that
$$
\int_0^t\phi(s)\,ds\le\int_0^T\|\nabla v\|^2\,dt+C(\Omega,T,f,g),
$$
with a bound uniform in $\delta$, for $\delta$ small. Moreover
\begin{align*}
\int_0^t\psi(s)\,ds
\le&
C\left(\int_0^T\|\nabla\Ud\|^2\right)^{\frac13}
\left(\int_0^T\|\nabla\zd\|^2\right)^{\frac23}\|\Ud\|_{L^\infty}\\  
   &+C\left(\int_0^T\|\nabla\Ud\|^2\right)^{\frac23}
\left(\int_0^T\|\nabla\zd\|^2\right)^{\frac13}\|\zd\|_{L^\infty} 
    +C\|\zd\|_{L^\infty}\int_0^T\|\nabla\zd\|^2\\
\le& C(\delta^{\frac23}+\delta^{\frac56}+\delta^{\frac32})\\
\le& C\delta^{\frac23},
\end{align*}
so that, by~\eqref{limest} and Gronwall's inequality, it follows that
\begin{align*}
\sup_{0\le t\le
  T}\|\od\|^2+\int_0^T\Big(\|\D(\od)\|^2+
\frac1\delta\|\od\|_\Gamma^2\Big)\,dt 
&\le\int_0^T[\delta\chi(s)+\psi(s)]\e^{\int_s^T\phi(r)\,dr}\,ds\\
&\le\e^{\int_0^T\phi(s)\,ds}\int_0^T[\delta\chi(s)+\psi(s))]\,ds\\
&\le C\delta^{\frac23}
\end{align*}
and finally
$$
\sup_{0\le t\le T}\|\od\|^2+
\int_0^T\Big(\|\D(\od)\|^2+\frac1\delta\|\od\|_\Gamma^2\Big)\,dt
=\mathcal{O}(\delta^{\frac23}). 
$$
Since the same norms of $\zd$ are of order $\delta^{\frac12}$, by the
previous lemma, we finally deduce that $u_\delta-v$ is of the order
$\delta^{\frac13}$, and the theorem is proved.
\begin{remark}
  {\rm The tangential trace of the function $\od$ converges a little
  better, since from the above estimate we can deduce
$$
\|\od\|_{L^2(\Gamma\times(0,T))}=\mathcal{O}(\delta^{5/6}).
$$}
\end{remark}
\begin{remark} {\rm The estimates show that there is a crucial
loss in the estimates due 1) to the boundary effect, in the
estimate~\eqref{stimaG}; 2) to the non-linearity (recall the
contribution of $\Ud$ in the estimate of $\psi(t)$). We observe that
the original model was supplemented by the following boundary
condition 
$$
u\cdot\nv=\delta^2 g,
$$
but we used in fact $u\cdot\nv=\delta\, g$. By running through the
proof of Theorem~\ref{mainlim} one gets, using the boundary condition
$$
u\cdot\nv=\delta^\alpha g,
$$
that $\zd=\mathcal{O}(\delta^{\alpha-\frac12})$, so that
$\int\psi(t)\,dt=\mathcal{O}(\delta^{\frac43(\delta-\frac12)})$ 
and finally $\od=\mathcal{O}(\delta^{\frac23(\alpha-\frac12)})$.
In particular, for $\alpha=2$, that is, the value corresponding
to the original model, the behaviour of $\zd$ matches the loss
in the estimates due to the boundary term~\eqref{boundary}. Hence,
for larger values of $\alpha$, the convergence remains slow, because
of this term. Supposedly, the loss in the convergence rate caused by
this term is an intrinsic feature of the problem, while those caused
by the non-linearity may have technical reasons.}
\end{remark} 
\bibliography{BerselliRomito}
\bibliographystyle{siam}
\end{document}